%% file: main.tex
\documentclass[a4paper, 11pt]{article}

\usepackage[english]{babel}
\usepackage[utf8]{inputenc}
\usepackage{subfiles}
\usepackage[a4paper,margin=2.5 cm]{geometry}
\usepackage{csquotes}
\usepackage[numbib,nottoc]{tocbibind}
\usepackage{amsmath}
\usepackage{indentfirst}
\usepackage{amssymb}
\usepackage{graphicx}
\usepackage{adjustbox}
\usepackage{tabularx}
\usepackage{booktabs}
\usepackage{tabu}
\usepackage{caption}
\usepackage{blindtext}
\usepackage[dvipsnames]{xcolor}
\usepackage[colorinlistoftodos]{todonotes}
\usepackage[colorlinks=true, allcolors=Cyan]{hyperref}
\usepackage{listings}
\usepackage{color}
\usepackage{setspace}
\usepackage{shapepar}
\usepackage{graphicx}
\usepackage{ragged2e}
\usepackage{listings}
\usepackage{algorithmicx}
\usepackage{algorithm}
\usepackage{float}
\usepackage{color,soul}
\usepackage{breakurl}
\usepackage{enumitem}

\usepackage[noend]{algpseudocode}
\usepackage{listings} \usepackage{color} 

\usepackage{amsfonts,amssymb,amsthm,epsfig,epstopdf,titling,url,array,tikz,graphicx,tkz-euclide,mathtools}
\newcommand{\cmnt}[1]{\ignorespaces}
\DeclarePairedDelimiter\ceil{\lceil}{\rceil}

\DeclareMathOperator{\spn}{span}

\makeatletter
\newcommand{\tpmod}[1]{{\@displayfalse\pmod{#1}}}
\makeatother

\makeatletter
\newtheorem*{rep@theorem}{\rep@title}
\newcommand{\newreptheorem}[2]{%
\newenvironment{rep#1}[1]{%
 \def\rep@title{#2 \ref{##1}}%
 \begin{rep@theorem}}%
 {\end{rep@theorem}}}
\makeatother

\theoremstyle{plain}
\newtheorem{theorem}{Theorem}[section]
\newreptheorem{theorem}{Theorem}
\newtheorem{lemma}{Lemma}[section]
\newreptheorem{lemma}{Lemma}
\newtheorem{prop}{Proposition}[section]

\theoremstyle{definition}
\newtheorem{definition}{Definition}[section]
\newreptheorem{definition}{Definition}

\newtheorem{remark}{Remark}[section]

\usepackage[backend=biber]{biblatex}
\title{Asymptotics of $d$-Dimensional Visibility}
\author{Ezra Erives \and Srinivasan Sathiamurthy \and Zarathustra Brady}
\date{August 2019}

\begin{document}

\maketitle

\subfile{Sections/Abstract}

\subfile{Sections/S1}
\subfile{Sections/S2}
\subfile{Sections/S3}
\subfile{Sections/S4}
\subfile{Sections/S5}
\subfile{Sections/S6}
\subfile{Sections/S7}
\subfile{Sections/S8}

\subfile{Sections/Future_Work}
\subfile{Sections/Acknowledgements}

\printbibliography
\end{document}

%% file: Sections/Abstract.tex
\begin{abstract}
    We consider the space $[0,n]^3$, imagined as a three dimensional, axis-aligned grid world partitioned into $n^3$ $1\times 1 \times 1$ unit cubes. Each cube is either considered to be empty, in which case a line of sight can pass through it, or obstructing, in which case no line of sight can pass through it. From a given position, some of these obstructing cubes block one’s view of other obstructing cubes, leading to the following extremal problem: What is the largest number of obstructing cubes that can be simultaneously visible from the surface of an observer cube, over all possible choices of which cubes of $[0,n]^3$ are obstructing? We construct an example of a configuration in which $\Omega\big(n^\frac{8}{3}\big)$ obstructing cubes are visible, and generalize this to an example with $\Omega\big(n^{d-\tfrac{1}{d}}\big)$ visible obstructing hypercubes for dimension $d>3$. Using Fourier analytic techniques, we prove an $O\big(n^{d-\tfrac{1}{d}}\log n\big)$ upper bound in a reduced visibility setting.
\end{abstract}

%% file: Sections/S1.tex
\section{Introduction}
\subsection{Visibility Problems}
Consider a configuration of (opaque) objects in space. Two objects are said to be visible from each other if there exists an unobstructed line segment between a point on the first object and a point on the second. We are curious about the maximum number of objects which may be visible from a particular point in space in a worst case scenario. We consider the simple case in which our objects are unit cubes with vertices at integer coordinates, bounded between $0$ and $n$. It is clear that there is a configuration in which you can see at least a quadratic (in $n$) number of cubes: you can see all of the cubes that share a face with the boundary of the $n\times n\times n$ cube simultaneously if there are no other obstructing cubes inside of the grid. There are also clearly a maximum of $n^3$ cubes that lie in your range of visibility. However, it is not clear if it is possible to see a number of obstructing cubes that is cubic in $n$.

Although questions of a similar flavor have been asked before, the techniques used to solve them are inapplicable to our context. One famous family of questions (see \cite{orchard_problem}), namely the Orchard Visibility Problem and its variants, considers a circular orchard bounded by a radius $R$ in which a tree of radius $r<1$ is centered at every lattice point. Both the observer and the center of the orchard are located at the origin. The problem asks what the maximum size of $r$ is for which there exists a line of sight connecting the observer to a point on the boundary of the orchard. In \cite{OVP}, Kruskal demonstrates that there exists such a line of sight if and only if $r < \frac{1}{R^2+1}$. 
 
 A more general view obstruction problem studied in \cite{general_obstruction}, in which centered at each point in the set $(-\tfrac{1}{2},-\tfrac{1}{2}\dots,-\tfrac{1}{2})+\mathbb{N}^d$ in $\mathbb{R}^n$ is a centrally symmetric convex body. The bodies are expanded uniformly until they block all rays emanating from the origin and into the open positive cone. The problem has been solved for balls in dimensions $d=2,3,4$.
 
 A related problem considers the observer to be positioned at the origin, and trees to be located at lattice points (with radius $0$). The maximum number of trees visible is in bijection with the number of pairs $(x,y)\in \mathbb{N}^2$ such that $\gcd(x,y)=1$, which for an infinitely large orchard is $\frac{6}{\pi^2}$ of the trees in the orchard. For $d>2$ dimensions, this generalizes to $\frac{1}{\zeta(d)}$ of the trees, where $\zeta$ is the Riemann zeta function.

\begin{figure}
    \centering
    \includegraphics[scale=0.35]{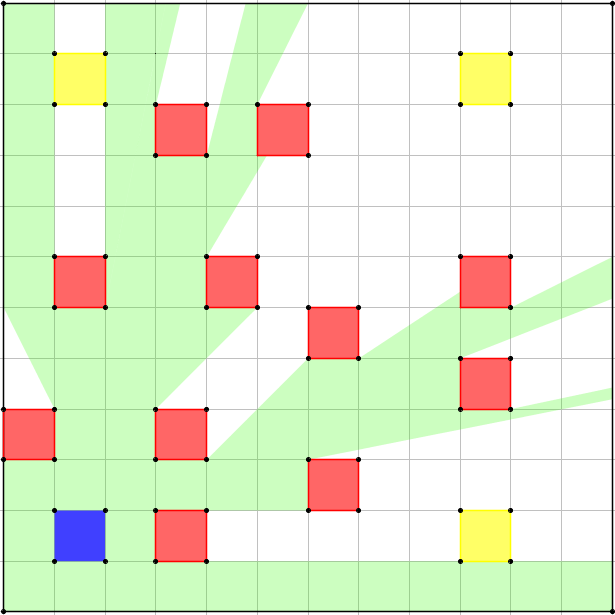}
    \caption{The locus of visible points from the darkened (blue) square in the lower left corner.}
    \label{fig1}
\end{figure}

The two dimensional version of our problem has been solved in \cite{zeb}, where Brady considers an $n \times n$ axis-aligned grid. Figure \ref{fig1} shows visibility from the darkened blue square in the lower left hand corner, with the obstructing squares visible from the blue square colored red, and the ones not visible colored yellow. The locus of points visible to the blue square are shaded in green.
     
Brady asks: If the number and placement of the obstructing squares in the grid is optimal, then what is the largest number of obstructing squares, as a function of $n$, that can be visible to a given square? 
Brady used elementary techniques to demonstrate that the answer to this problem is $\Theta(n\sqrt n)$. To do so, he split the $n \times n$ grid into (mostly) disjoint parallelograms, and computed lower and upper bounds on this maximal value. While the elementary approach used there doesn't generalize to higher dimensions, we do use the same parallelogram approach in our argument.
     
\begin{figure}[H]
        \centering
        \begin{minipage}{.48\textwidth}
        \centering
        \includegraphics[width=.7\linewidth]{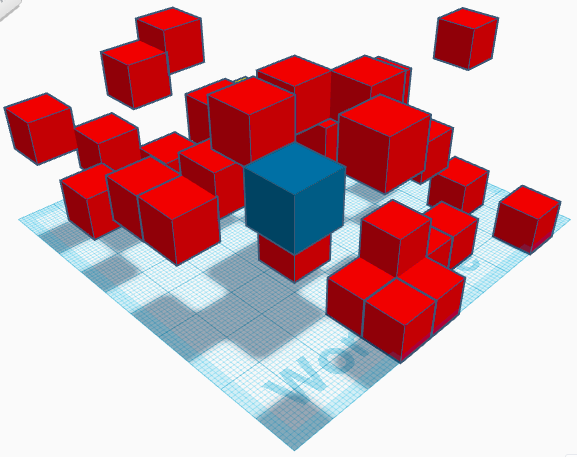}
        \end{minipage}
        \begin{minipage}{.48\textwidth}
        \centering
        \includegraphics[width=.7\linewidth]{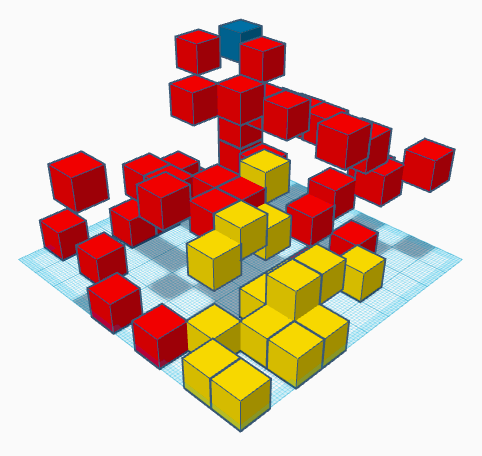}
        \end{minipage}
    \caption{Visibility is taken from perspective of the dark blue cube. The cubes that are both obstructing and visible to the blue cube are painted red while the non-visible obstructing cubes are painted yellow.}
    \label{fig2}
\end{figure}

In our paper, we generalize the two dimensional bounds Brady obtained to $d>2$ dimensions. In doing so, the two dimensional $n \times n$ grid of squares becomes a $d$-hypercube of side length $n$ consisting of $n^d$ unit hypercubes, each of which is either empty or obstructing. Within this larger $d$-hypercube of side length $n$, we seek the maximum possible number of obstructing unit $d$-hypercubes visible from a given obstructing unit $d$-hypercube. It is easiest to visualize this question when $d=3$, and so we shall go about analyzing the problem in three dimensions before extending our results to higher numbers of dimensions. Figure \ref{fig2} illustrates visibility from the dark blue cube in the case of $d=3$.
    
\subsection{Main Results}

\cmnt{
3.1,
4.1,
4.2,
7.1,
7.2,
8.1,
}

Adopting a similar argument to that taken in the two-dimensional case (see \cite{zeb}), we assume the observer is the cube adjacent to the origin  and divide $[0,n]^3$ into $1\times 1\times n$ parallelepipeds through the origin to construct a lower bound. By projecting the possible obstructing cubes intersecting the parallelepipeds' long edges onto the bottom faces of the parallelepipeds, we construct a partially ordered set that characterizes the conditions under which obstructions can block each other. The task of constructing sets of simultaneously visible cubes from the origin is then transformed into one of constructing an antichain of maximal size of a certain partially ordered set. We demonstrate the existence of ``small'' vectors modulo $n$ that span a lattice corresponding to an antichain of the partially ordered set.

\begin{reptheorem}{thm31}
There is a configuration in which the number of obstructing unit cubes within a cube of side length $n$ visible from the origin is at least $\Omega(n^\frac{8}{3})$.
\end{reptheorem}

The partially ordered set described above is then modified to model visibility in $d>3$ dimensions. We extend the techniques used in three dimensions to construct a large antichain (coming from a $d-1$ dimensional lattice) of the generalized partially ordered set.
 
\begin{reptheorem}{thm41.5}\label{repthm41.5}
Let $\vec{t} = (t_1, ..., t_{d-1}) \in \{1,...,p-1\}^{d-1}$. Let $S_{\vec{t}}$ denote the set of $2^{d-1}$ partially ordered sets, each of which is of the form
\[
\big\{ \big( (\pm t_1 \cdot k) \tpmod p, \dots, (\pm t_{d-1} \cdot k)\tpmod p, k
\big)\mid 0\le k <p\big\},
\]
under product order, for some fixed choice of signs.
Then for each $\vec{t}$ there exists an element of $S_{\vec{t}}$ whose width is $\Omega (p^{1-\frac{1}{d}})$.
\end{reptheorem}


The linear (modulo $p$) structure of the posets considered in Theorem 4.1 is crucial for dimension $d > 2$. For any $d$-dimensional poset, one can define a similar family of $2^{d-1}$ related posets by reversing the order on some of the coordinates, and one might be tempted to believe that an analogue of Theorem \ref{thm41.5} holds for any such situation. When $d = 2$, this is indeed true, and it follows easily from the well known Erdős–Szekeres theorem. However, in dimension $d > 2$ this is no longer the case. In \cite{erdos/szekeres}, Szabó and Tardos consider multidimensional generalizations of the Erdős–Szekeres theorem and demonstrate the existence of situations in which all $2^{d-1}$ partially ordered sets in the family of posets have width $O(n^{e_d})$ with $e_d < 1-\tfrac{1}{d}$ for $d \ge 3$ (for instance, their $e_3$ is $\tfrac{5}{8}$, which is less than $\tfrac{2}{3}$).

The existence of the partially ordered set from Theorem \ref{thm41.5} leads to a construction of a configuration in which many obstructions are visible.
 
\begin{reptheorem}{thm42}
There is a configuration in which the number of obstructing unit $d$-hypercubes withint $[0,n]^d$ which are visible from the origin is at least $\Omega(n^{d-\tfrac{1}{d}})$.
\end{reptheorem}

As an aside, we note that the $\Omega(n^{1-\tfrac{1}{d}})$ lower bound achieved above matches with the width of a \emph{random} $d$-dimensional partial order of size $n$, as computed by Brightwell \cite{brightwell} (see discussion of Theorem \ref{thm22} for more).

In approaching an upper bound on the number of visible obstructions, we first consider a reduced visibility environment, in which we restrict visibility to only lines of sight parallel to the edges of the $d$-parallelepiped in consideration. Using the same partially ordered set as used in the lower bound in three dimensions, we demonstrate that there exists a chain cover of sufficiently small size. We do so by studying the value $h_p$ defined as follows.

\begin{repdefinition}{def51}  For prime $p$, positive integer $d>2$ and $\vec t=(t_1,\dots,t_{d-1},1)\in \mathbb{Z}_p^d$, we define:
\[h_{p}(\vec{t}) := \min_{0 \le a < p} \max ((at_1)\tpmod p,\dots,(at_{d-1})\tpmod p, a).\]
\end{repdefinition}


 \begin{reptheorem}{thm74}
 The average value of $h_p(\vec{t})$ as $\vec{t}$ varies is bounded by $O(p^{\tfrac{d-1}{d}} \log p)$, that is,
 \[
\mathbb{E}_{\vec{t}}[h_p(\vec{t})] \ll p^{\tfrac{d-1}{d}} \log p.
\]
 \end{reptheorem}
 
 
 \begin{reptheorem}{thm71}
 The largest number of cubes visible in the $d$ dimensional toy upper bound visibility environment is $O(p^{d-\tfrac{1}{d}}\log p)$. 
 \end{reptheorem}
 
In order to prove these results, we introduce a dual height $h_p^*$ which is small when there is a simple reason for $h_p$ to be large.

\begin{repdefinition}{dualhp}For $\vec{t} \in \mathbb{Z}_p^d$, we define the \emph{dual height} $h_p^*(\vec{t})$ by
\[
h_p^*(\vec{t}) = \displaystyle{\min_{\substack{\vec \alpha\cdot \vec{t} \equiv 0 \\ \alpha \not \equiv \vec 0}} \max |\alpha_i|}.
\]
\end{repdefinition}

Note that if the $t_i$ satisfy a simple linear relation such as $at_1 + bt_2 + ct_3 = 0 \tpmod{p}$ with $a,b,c$ small positive constants, then $h_p(t)$ must be at least $\frac{p}{(a+b+c)} \sim \frac{p}{max(a,b,c)}$ (up to a factor of $d$). We prove a weak converse to this.

\begin{replemma}{lemhp} For all $\vec{t} \in \mathbb{Z}_p^d$, we have
\[
h_p(\vec{t}) \ll_d \frac{p\log p}{h_p^*(\vec{t})}.
\]
\end{replemma}
 
 Up to this point, lines of sight under consideration were restricted to just those parallel to the lateral edges of the current parallelotope. We next examine visibility in an environment where this restriction is no longer in place. As we were not able to solve for an upper bound in an unrestricted visibility environment, we weaken the problem.
 
 \begin{repdefinition}{def81}
  We say that a $d$-dimensional hypercube blocks a ray of light \emph{at the angle} $\theta$ if the ray of light intersects some $d-1$-dimensional facet of the hypercube at an angle at most $90^{\circ}-\theta$ away from the normal vector to that facet.
 \end{repdefinition}
 
 \begin{reptheorem}{thm81}
 The largest number of visible obstructions, in the setting where light fails to interact with any obstruction that does not block it at the angle $45^\circ$, is at most $O(n^{d-\tfrac{1}{d}}\log n)$.
 \end{reptheorem}

 \subsection{Organization of Material}
 In Section \ref{sec2}, a brief introduction is given to partially ordered sets. In Section \ref{sec3}, a lower bound on the number of obstructing cubes visible in three dimensions is proven. In Section \ref{sec4}, the bound presented in Section \ref{sec3} is generalized to $d>3$ dimensions. Section \ref{sec5}  provides an introduction to the so called toy upper bound, a simplification of the true upper bound. In Section \ref{sec6}, a brief introduction is given to the discrete Fourier transform. In Section \ref{sec7}, the results of Sections \ref{sec5} and \ref{sec6} are combined to present a bound on visibility in the restricted setting of the toy upper bound in $d>2$ dimensions. Finally in Section \ref{sec8}, an upper bound is presented in a shallow light visibility environment which is stronger than the setting of Section \ref{sec7} but weaker than the full visibility setting we are interested in.
 
 
 

%% file: Sections/S2.tex
\section{Partially Ordered Sets} \label{sec2}


\begin{definition}
Let $P$ be a poset with relation $\le$. A $\emph{chain}$ of $P$ is a subset $S\subseteq P$ such that for all $a,b\in S$, either $a\le b$ or $b \le a$.
\end{definition}
\begin{definition}
Let $P$ be a poset with relation $\le$. A $\emph{antichain}$ of $P$ is a subset $S\subseteq P$ such that for all $a\ne b\in S$, neither $a\le b$ nor $b\le a$.
\end{definition}
\begin{definition}
The \emph{width} of a finite poset is the size of its largest antichain.
\end{definition}
\begin{theorem}[Dilworth's Theorem \cite{dilworth}]
The width of a finite poset $P$ is equal to the minimum number of chains into which $P$ can be partitioned. 
\end{theorem}


We say that a total ordering $\le_1$ on a poset $P$ is \emph{compatible} with a partial ordering $\le_2$ if for all $a,b\in P$, $a\le_2 b$ implies $a\le_1 b$.
\begin{definition}
A \emph{linear extension} of a poset $P$ is a total ordering of $P$ which is compatible with the partial order on $P$.
\end{definition}
\begin{definition}
The dimension of a poset $P$ with partial order $\le$ is the least integer $d$ for which there exists a family $R=(\le_1,\le_2,\dots,\le_d)$ of linear extensions of $P$ such that 
\[\le\ =\bigcap_{i=1}^d \le_i,\]
where $\le, \le_i$ are treated as subsets of $P\times P$.
\end{definition}
\begin{prop}
The partially ordered set $S_d=\{(x_1,x_2,\dots,x_d)\mid x_i\in \mathbb{R}\}$ under product order has dimension $d$.
\end{prop}

As a result of the above proposition, we refer to posets whose elements are $d$-tuples as being $d$-dimensional.

As we will see, the question of visibility reduces to one of analyzing certain partially ordered sets. We conclude this section with a result due to Brightwell \cite{brightwell}, which provides upper and lower bounds on the width of a \emph{random} product-ordered $d$-dimensional tuple, which the reader can compare to the bounds we will prove later on the widths of the specific posets we are interested in.

\begin{theorem}[Brightwell] \label{thm22}
There exists a constant $C$ such that, for each fixed $d$, almost every $P_d(n)$ satisfies
\[
\Bigg( \frac{1}{2}\sqrt{d}-C \Bigg)n^{1-\frac{1}{d}} \leq W_d(n) \leq \frac{7}{2}dn^{1-\frac{1}{d}},
\]
where $P_d(n)$ denotes a random $d$-dimensional poset under product order, and where $W_d(n)$ is the width of such a poset.
\end{theorem}



%% file: Sections/S3.tex
\section{A Lower Bound Construction in Three Dimensions}\label{sec3}
In this section, we prove the existence of a set of $\Omega(n^{\tfrac{8}{3}})$ obstructing cubes, all of which are simultaneously visible from an observer cube centered at the origin. In order to simplify some of the number theoretic computation, we will assume that $n$ is some prime $p$ by possibly replacing $n$ with the largest prime $p$ which is less than $n$. Bertrand's postulate ensures that in doing so, our bound remains unchanged asymptotically. 


\subsection{Setup}\label{sec31}

\begin{figure}
    \centering
    \begin{minipage}{.48\textwidth}
    \centering
    \includegraphics[width=.8\linewidth]{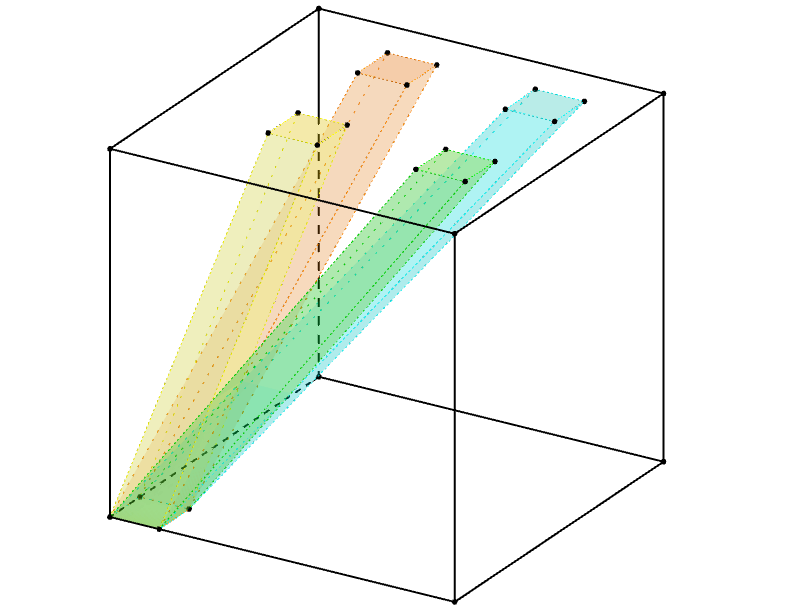}
    \caption{Example parallelepipeds}
    \label{fig3}
    \end{minipage}
    \begin{minipage}{.48\textwidth}
    \centering
    \includegraphics[width=.8\linewidth]{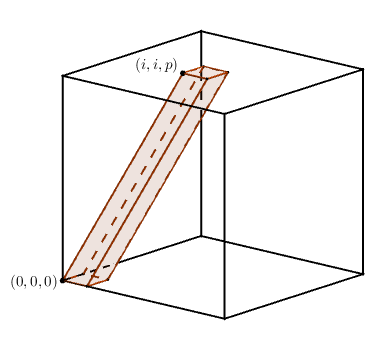}
    \caption{The parallelepiped $\mathcal{P}_{i,j}$ with characteristic vertex $(i,j,p).$}
    \label{fig4}
    \end{minipage}
\end{figure}
    
We first consider the set of parallelepipeds with opposite and parallel square faces, one of which is a unit square whose vertices have integer coordinates on the upper face of the cube, and the other the unit square on the bottom face of the cube with one vertex at the origin. There are $p^2$ such parallelepipeds, one for each unit square on the cube's upper face. Specifically, we are considering parallelepipeds with vertices at coordinates $(0,0,0)$, $(0,1,0)$, $(1,0,0)$, $(1,1,0)$, $(i,j,p)$, $(i+1,j,p)$, $(i,j+1,p)$ and $(i+1,j+1,p)$ where $p$ is the size of the grid, $p$ is a prime number, and $0 \leq i, j \leq p-1$. 

We now shift our focus to one of these parallelepipeds, call it $\mathcal{P}_{i,j}$, which we shall refer to as simply $\mathcal{P}$ when there is no danger of confusion, whose top face has vertex $(i,j,p)$ closest to the origin. We will refer to $(i,j,p)$ as the $\textit{characteristic vertex}$ of $\mathcal{P}$. Let $E_{\mathcal{P}_{i,j}}$, or more simply $E_\mathcal{P}$ be the edge of $\mathcal{P}_{i,j}$ containing both $(i,j,p)$ and the origin. Denote $E_\mathcal{P}$, as well as the other three edges of $\mathcal{P}$ parallel to $E_\mathcal{P}$ as the $\textit{lateral edges}$ of $\mathcal{P}$. 

Within $\mathcal{P}$, we will only consider the set of possible obstructing cubes intersecting $E_\mathcal{P}$. We additionally restrict visibility to lines of sight parallel to the lateral edges of $\mathcal{P}$. As a result, a line of sight passing through the observer cube is equivalent to the line of sight passing through the bottom face of the observer cube, so we may flatten the observer cube to merely its bottom face without sacrificing any visibility.

The crucial observation is that from the perspective of the observer face, all possible obstructing cubes intersecting the line $E_\mathcal{P}$ can be seen only by their bottom face (note that we are not considering all cubes intersecting $\mathcal{P}$, just those which intersect $E_\mathcal{P}$ - for the cubes intersecting $\mathcal{P}$ but not $E_\mathcal{P}$, this crucial observation is not true). It follows that a set of obstructing cubes, all of which intersect $E_\mathcal{P}$, are simultaneously visible to the observer's face if and only if the same statement is true of their bottom faces. 

Note that obstructing cubes near the origin have the potential to be counted as intersecting the lateral edge containing the characteristic vertex for more than some fixed constant number of times as we enumerate over all such $\mathcal{P}$. For this reason, we will further restrict candidates for obstructing cubes to those intersecting the edge $E_\mathcal{P}$ of $\mathcal{P}$ that are in the upper half of the parallelepiped, ensuring that each obstructing cube is counted in only a fixed constant number of parallelepipeds. In the following argument however, we will assume that \textit{any} obstructing cube along $E_\mathcal{P}$ is a candidate for being obstructing. This simplification is justified by Proposition \ref{prop31}, stating that the width of the partially ordered set modelling visibility for \textit{all} cubes along $E_\mathcal{P}$ is at most a constant factor larger than the partially ordered set modelling visibility along just the upper half of $E_\mathcal{P}$.

\subsection{A Partially Ordered Set for Visibility Along $E_\mathcal{P}$}\label{sec32}
In this subsection we construct a partially ordered set for $\mathcal{P}_{i,j}$ that models visibility along the edge $E_\mathcal{P}$, and whose width is precisely the largest possible number of simultaneously visible obstructions along $E_\mathcal{P}$. By switching our choice of $i$ and $j$, the construction can be extended to any such parallelepiped of the form described in Subsection $\ref{sec31}$. In particular, we will show that among the four partially ordered sets attached $\mathcal{P}$ and the parallelepipeds with characteristic vertices $(p-i,j,p)$, $(i,p-j,p)$ and $(p-i,p-j,p)$ respectively, at least one has sufficiently large width.


Consider two square faces at heights $k_1$ and $k_2$ intersecting $E_\mathcal{P}$ and the origin, respectively. When these two squares are projected onto the observer face, the corners are taken to the points $(1- \{ \tfrac{i k_1}{p} \},1- \{ \tfrac{j  k_1}{p} \}, 0 )$ and $(1- \{ \tfrac{i  k_2}{p} \},1- \{ \tfrac{j  k_2}{p} \}, 0 )$, respectively. It follows that the face at height $k_1$ is visible from the observer face if either one of the following conditions are true.
\begin{enumerate}
    \item $k_1<k_2$ (face $k_1$ is lower than face $k_2$),
    \item $1-\{\tfrac{i k_1}{p} \} > 1-\{ \tfrac{i  k_2}{p} \}$  (the corner face $k_1$ ``sticks out" from behind face $k_2$ with respect to $x$ coordinate), 
    \item $1-\{\tfrac{j k_1}{p}\} > 1-\{\tfrac{j  k_2}{p}\}$ (the corner face $k_1$ ``sticks out" from behind face $k_2$ with respect to $y$ coordinate).
\end{enumerate}
Thus, the face at height $k_2$ can only block the observer face's view of the face at height $k_1$ if every coordinate of $( 1- \{ \tfrac{i k_2}{p} \},1-\{ \tfrac{j k_2}{p} \}, p-k_2)$ is larger than the corresponding coordinate in $(1- \{ \tfrac{i  k_1}{p} \},1- \{ \tfrac{j  k_1}{p} \}, p-k_1)$. As a result, a set of obstructing faces all intersecting $E_\mathcal{P}$, all of whose elements are simultaneously visible from the observer face, corresponds to an antichain of the partially ordered set $\big\{(\{ \tfrac{i k}{p} \}, \{ \tfrac{j  k}{p} \}, k ) \mid 0\le k < p \big\}$, under product order. This poset is isomorphic to the much simpler $(ik\pmod p, jk \pmod p, k) \mid 0\le k < p)$ which we shall refer to as $S_{i,j}$. The width of this poset is the maximum number of obstructing cubes intersecting $E_\mathcal{P}$ that are all simultaneously visible from the observing face with respect to the restricted lines of sight. 

Now consider $\mathcal{P}_{p-i,j}$, $\mathcal{P}_{i,p-j}$ and $\mathcal{P}_{p-i,p-j}$. These parallelepipeds have posets analogous to $S_{i,j}$, with elements of the form $(-i k\pmod p, jk \pmod p, k)$, $(i k\pmod p, -jk \pmod p, k)$, and $(-i k\pmod p, -jk \pmod p, k)$, and which we will denote as $S_{p-i,j}$, $S_{i,p-j}$ and $S_{p-i,p-j}$ respectively. It is then only natural to consider these three partially ordered sets together with $S_{i,j}$ and alternatively their four corresponding parallelepipeds as belonging to the same \textit{family}. In this manner, the entire set of such parallelepipeds, the size of which is quadratic in $p$, may be partitioned into families of four parallelepipeds.

\begin{prop}\label{prop31}
Let $p>2$ be a prime, $(i,j)\in \mathbb{Z}_p^2$, and let $S$ be one of the four partially ordered sets of the form $\{(\pm (ki)\pmod p, \pm (kj) \pmod p, k)\mid 0 \le k < p\}$ for some fixed choice of signs, and let $w$ be the width of such a poset. Denote $S_{-}$ and $S_{+}$ to be the subsets of $S$ for which $\frac{p-1}{2}<k<p$ and $0\le k < \frac{p-1}{2}$ respectively, and let $w_{-}$ and $w_{+}$ the the two posets' respective widths. Then $w_-=w_+\ge\tfrac{w-1}{2}$.
\end{prop}


    The proof of Proposition \ref{prop31} can easily be extended to hold true for any partially ordered set of the form $\{(\pm t_1k\pmod p,\dots, \pm t_m k \pmod p, k)\mid 0 \le k < p\}$ for some fixed choice of signs and taken under product order, where $m$ any positive integer.

\begin{prop}\label{prop32}
Let $p>2$ be a prime, $(t_1,\cdots,t_{d-1})\in \mathbb{Z}_p^{d-1}$, and let $S$ be one of the $2^{d-1}$ partially ordered sets of the form $\{(\pm (kt_1)\pmod p, \cdots, \pm (kt_{d-1}) \pmod p, k)\mid 0 \le k < p\}$ for some fixed choice of signs, and let $w$ be the width of such a poset. Denote $S_{-}$ and $S_{+}$ to be the subsets of $S$ for which $\frac{p-1}{2}<k<p$ and $0\le k < \frac{p-1}{2}$ respectively, and let $w_{+}$ and $w_{-}$ the the two posets' respective widths. Then $w_-=w_+\ge\tfrac{w-1}{2}$.
\end{prop}

\begin{figure}
    \centering
    \includegraphics[scale=0.2]{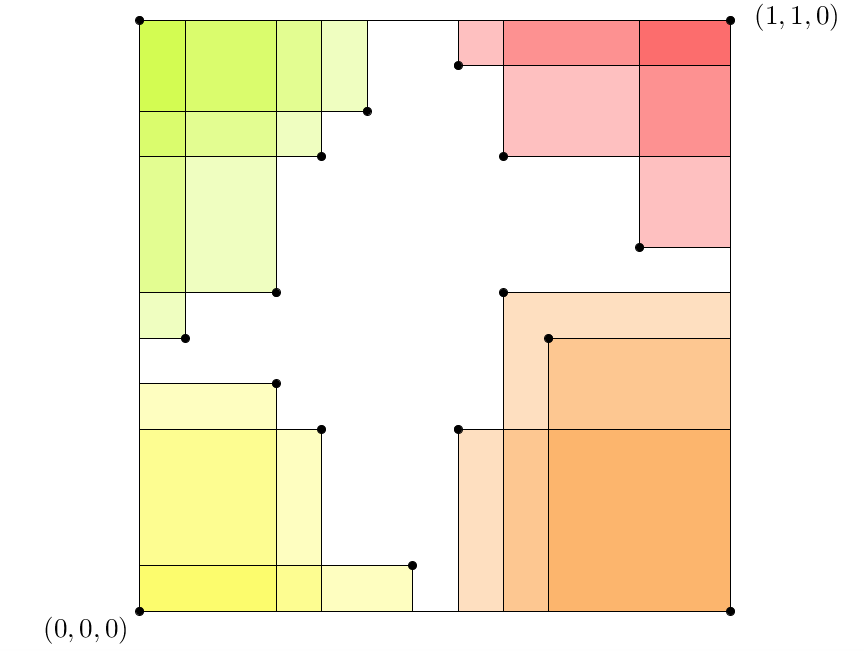}
    \caption{Bottom faces of obstructing cubes projected onto the observer cube and colored according to the lateral edge they intersect. The images of projections along $E_\mathcal{P}$ are shown in yellow.}
    \label{fig6}
\end{figure}

\subsection{Bounding the Width}
 Our objective now is to find a lower bound on the maximum width among the four posets within the family of $\mathcal{P}$. To bound the width, we show that there exists an antichain of one of the four previously mentioned partially ordered sets of sufficiently large size, and do so in a manner motivated by the following observation. Each element of $S_{i,j}$ may be viewed as a point within $[0,p)^3$. Note that just as these points in space are elements of our partially ordered sets, so too are the vectors obtained by taking the difference between any two of these points. From this perspective, we may consider a (shifted) two dimensional lattice of points from the poset as being generated by one element of our poset viewed as a starting point, and two more viewed as vectors. It follows that the intersection of $[0,p)^3$ and any such lattice of points from $S_{i,j}$ whose normal vector is of uniform sign corresponds to a maximal antichain of $S_{i,j}$, as the difference between any two elements of this plane is a vector of mixed sign.
 
 As we are working not just with $S_{i,j}$ but with $S_{p-i,j}$, $S_{i,p-j}$, and $S_{p-i,p-j}$ as well it suffices to merely find a lattice (of rank two) within one of the four posets of sufficiently large size. As any lattice within one of the four partially ordered sets exists within the other three albeit with different signs, it remains after such a lattice is found to simply choose the poset whose signs will guarantee the plane to have a normal vector all of whose coordinates are of the same sign.
 
In the spirit of this, we construct two linearly independent vectors in $S_{i,j}$ with all coordinates as small as possible absolute value wise, and use the lattice spanned by the two vectors to make a statement about the width of one of the four partially ordered sets.
 
 \begin{lemma}\label{lem32}
There exists a $v \in S_{i,j}$ for which the absolute value of all of $v$'s coordinates are $O\big(p^{\frac{2}{3}}\big)$.
\end{lemma}
\begin{proof}
We use a pigeonhole argument. Imagine each element of $S_{i,j}$ as a point within a cube of side length $p$. Divide the cube into an $m\times m \times m$ grid of cubes, with each cube having side length $\tfrac{p}{m}$. If $m$ satisfies $m^3<p$, then at least one of these cubes must have two points inside it. Take one of these cubes with two points inside of it, call the two points $v_1$ and $v_2$. Taking $v_3=v_1-v_2$, we see that all of $v_3$'s coordinates are at most the side length of the box. 
Letting $m=\big \lfloor p^{\frac{1}{3}} \big \rfloor > p^{\frac{1}{3}}-1$, gives a side length of $\tfrac{p}{m}=\tfrac{p}{\big\lfloor p^{\frac{1}{3}} \big\rfloor} < \tfrac{p}{p^{\frac{1}{3}}-1}=O\big(p^{\frac{2}{3}}\big)$.
\end{proof} 

We now construct a second vector in $S_{i,j}$ that is linearly independent from the one found in Lemma \ref{lem32} and that is of sufficiently small size.

\begin{lemma} \label{lem33}
There exist distinct and linearly independent $v_1,v_2\in S_{i,j}$, with largest coordinates $s_1$ and $s_2$ respectively, for which $s_1s_2= O\big(p^{\frac{4}{3}}\big)$.
\end{lemma}
\begin{proof}
 We first divide the $p \times p \times p$ cube into an $m\times m \times m$ grid of cubes, with each cube having side length $\tfrac{p}{m}$. By forcing $km^3<p$ for positive integer $k$, it is guaranteed that there is a cube $C$ with at least $k+1$ points from $S_{i,j}$ inside it. Simplifying, we have that $m<\big(\frac{p}{k}\big)^{\frac{1}{3}}$, so we take $m=\big \lfloor \big(\frac{p}{k}\big)^{\frac{1}{3}} \big \rfloor$ so that $\big(\frac{p}{k}\big)^{\frac{1}{3}}>m>\big(\frac{p}{k}\big)^{\frac{1}{3}}-1$.

 In finding a second vector, we must ensure that any new vector we find is not just a constant multiple of the first. To do so, it suffices to choose our $k$ so that $k+1$ is greater than the longest possible arithmetic sequence of points in $S_{i,j}$ that could be contained within $C$. Note that the side length of $C$ is $\frac{p}{m}\approx k^{\frac{1}{3}}p^{\frac{2}{3}}$. Let $s_1$ be the maximum coordinate of the vector $v_1$ in $S$ with the smallest maximum coordinate, absolute value wise (i.e. the worst possible case for generating an arithmetic sequence of points within $C$). Let $s_2$ be the maximum coordinate of the vector $v_2$ in $S$ which is not a multiple of $v_1$, which has the smallest maximum coordinate absolute value wise subject to not being a multiple of $v_1$. If $k$ satisfies $k > \frac{p}{ms_1} \approx \frac{k^{\frac{1}{3}}p^{\frac{2}{3}}}{s_1}$, then some pair of points in $C$ will have a difference which is not a multiple of $v_1$, so we will have $s_2$ at most the side length of $C$. We can take $k$ just a bit bigger than $ps_1^{-\frac{3}{2}}$. This implies that the side length of $C$ is about  $k^{\frac{1}{3}}p^{\frac{2}{3}}\approx \frac{p}{\sqrt{s_1}}$ and so $s_2= O\big(\frac{p}{\sqrt{s_1}}\big)$. By Lemma \ref{lem32}, $s_1=O\big(p^{\frac{2}{3}}\big)$ and so $s_1s_2=O(p\sqrt{s_1})=O\big(p^{\frac{4}{3}}\big)$, as desired.
 \end{proof}


 
 
 
 
 
 
 \begin{lemma}\label{lem34}
 When the elements of any one of $S_{i,j}$, $S_{p-i,j}$, $S_{i,p-j}$, and $S_{p-i,p-j}$ are taken as points in three dimensional space there exists a rank two lattice $\mathcal{L}$ of points from the chosen partially ordered set whose intersection with $[0,p)^3$ is of size $\Omega\big(p^{\frac{2}{3}}\big)$.
 \end{lemma}
 \begin{proof}
    We show that the result holds for $S_{i,j}$ which in turn proves the statement for all four posets by generality of lemmas $\ref{lem32}$ and $\ref{lem33}$.
    
    First observe that by Lemma $\ref{lem33}$, there exist $v_1,v_2\in S_{i,j}$ with largest coordinates $s_1$ and $s_2$ respectively for which $s_1s_2=O\big(p^{\tfrac{4}{3}}\big)$. To construct the desired lattice, we must now show that there is a point $\ell\in S_{i,j}$ suitably close to the center of $[0,p)^3$, after which we will take $\mathcal{L}=\{\ell + v \mid v\in \spn(v_1,v_2)\}$. Consider the union of intervals $I=(0,\frac{p}{6})\cup (\frac{5p}{6},p)$. Given some $v\in S_{i,j}$, the probability that any one of $v$'s three coordinates lies in $I$ is less than $\frac 1 3$. The probability that any of $v$'s three coordinates lie in $I$ is then strictly less than one (by the union bound) and so there must always exist an element of $S_{i,j}$ none of whose coordinates lie in $I$. We take $\ell$ to be this point none of whose coordinates lie in $I$, completing our construction of $\mathcal{L}$. The number of points inside $\mathcal{L}\cap [0,p)^3$ is up to a constant factor $\frac{p^2}{s_1s_2}=\Omega(p^\frac{2}{3})$, as desired.
\end{proof}


\begin{lemma}\label{lem35}
    At least one of $S_{i,j}$, $S_{p-i,j}$, $S_{i,p-j}$, and $S_{p-i,p-j}$ has width of $\Omega\big(p^{\tfrac{2}{3}}\big)$.
\end{lemma}
\begin{proof}
    Consider $S_{i,j}$. By Lemma \ref{lem34}, there exists a rank two lattice $\mathcal{L}_{i,j}$ whose intersection with $[0,p)^3$ contains only points from $S_{i,j}$ and is of size $\Omega(p^\frac{2}{3})$. The set of points $\mathcal{L}_{i,j}\cap[0,p)^3$ corresponds to an antichain if the normal vector to $\mathcal{L}_{i,j}$ has all coordinates strictly positive or strictly negative. The normal vector to $\mathcal{L}_{i,j}$ does not necessarily satisfy this sign requirement.
    
    Observe however that by negating (modulo $p$) the first coordinate of every point in the intersection $\mathcal{L}_{i,j}\cap[0,p)^3$, that is replacing the first coordinate $k$ of each point with $p-k$, we obtain a set of points all within $S_{p-i,j}$ that is precisely the intersection of a lattice $L_{p-i,j}$ with $[0,p)^3$. In particular, this new lattice has a corresponding normal vector whose sign on the first coordinate is the negative of that of the normal vector to $L_{i,j}$. 
    
    In this manner we may select one of the four partially ordered sets in the family of $S_{i,j}$ to force the signs of the coordinates of the normal vector to be uniform, guaranteeing that the final lattice is also an antichain. The result then follows from the fact that the lattice is of size $\Omega(p^\frac{2}{3})$.
\end{proof}

\begin{remark}
If this partially ordered set were to behave in the same manner as that of Theorem \ref{thm22}, its width would also be $\Omega(p^\frac{2}{3})$.
\end{remark}


\begin{prop}
For parallelepipeds $P_{i,j}$ and $P_{i',j'}$, if $|i-i'|\ge 6$ or $|j-j'|\ge 6$, then there exists no cube $C$ which intersects both $P_{i,j}$ and $P_{i',j'}$ and has $z$ coordinate at least $\frac{p}{2}$.
\end{prop}

\begin{theorem}\label{thm31}
There is a configuration in which the number of obstructing unit cubes within a cube of side length $n$ visible from the origin is at least $\Omega(n^\frac{8}{3})$.
\end{theorem}

\begin{proof}
By Bertrand's postulate, we may without loss of generality assume that $n=p$ is a prime. There are $\Theta(p^2)$ families of four parallelepipeds. By Lemma \ref{lem35} each such family has at least one element, call it $\mathcal{P}$, whose associated partially ordered set has width $\Omega(p^{\frac{2}{3}})$. As a result of Proposition \ref{prop31}, there is an antichain of this poset consisting only of points with third coordinate larger than $\tfrac{p-1}{2}$ also of size $\Omega(p^{\frac{2}{3}})$. By the bijection between sets of simultaneously visible obstructing cubes along the front lateral edge $E_\mathcal{P}$ of $\mathcal{P}$ and antichains of $\mathcal{P}$'s associated partially ordered set, there must exist a corresponding set of $\Omega(p^{\frac{2}{3}})$ simultaneously visible obstructing cubes along $E_\mathcal{P}$ all of which are in the upper half of the parallelepiped. Finally, enumerating over a collection of $\Theta(p^2)$ families of parallelepipeds whose characteristic vertices differ from each other by at least six in at least one coordinate yields a total of $\Omega(p^\frac{8}{3})$ total visible obstructing cubes.
\end{proof}



%% file: Sections/S4.tex
\section{A Lower Bound in $d>3$ Dimensions} \label{sec4}
\subsection{Generalization of the Three Dimensional Geometric Setup\cmnt{and a Representative Partially Ordered Set}}\label{sec41}
In this section, we generalize the results from Section $\ref{sec3}$ to $d>3$ dimensions. While in three dimensions, we asked how many obstructing cubes could possibly be seen from an observer cube, we now ask how many $d$-\emph{hypercubes}, which we will abbreviate to just hypercubes when there is no danger of ambiguity, can be seen from an observer $d$-\emph{hypercube}. The geometric approach is essentially the same.

As in Section $\ref{sec3}$, we assume that $n=p$ is a prime. We then consider the $d$-hypercube $C$ with side length $p$ formed by $[0,p)^d$ and assume that the observer is the unit $d$-hypercube adjacent to the origin.

Consider the vertex $(t_1,t_2,\dots,t_{d-1},p)$ on the upper $d-1$-dimensional facet of $C$. This is the $d$-dimensional analogue of the point $(i,j,p)$ (see Section \ref{sec3}). We may then construct the unique $d$-parallelotope $\mathcal{P}$ whose lower base is the bottom $d-1$-dimensional facet of the observer, and which contains the edge with endpoints at the origin and $(t_1,t_2,\dots,t_{d-1},p)$. As in Section \ref{sec3}, we refer to the point $(t_1,t_2,\dots,t_{d-1},p)$ as the characteristic vertex of $\mathcal{P}$. Additionally we will refer to the segment connecting the characteristic vertex of $\mathcal{P}$ to the origin as $E_\mathcal{P}$. For each such choice of $\mathcal{P}$, we will try and maximize the number of obstructing hypercubes intersecting $E_\mathcal{P}$ that are simultaneously visible from the observer. 

The edge $E_\mathcal{P}$ has the special property that any hypercube obstructing $E_\mathcal{P}$ is visible to the observer only by its bottom facet, and as a result, any set of obstructing hypercubes intersecting $E_\mathcal{P}$ is visible if and only if the corresponding set of bottom facets are also all simultaneously visible. In a near identical manner to Subsection $\ref{sec32}$, it can be seen that the largest number of simultaneously visible obstructing hypercubes along $E_\mathcal{P}$ is the width of the set 
\begin{equation*}\label{dag2}
\{((kt_1) \tpmod p,\dots, (kt_{d-1}) \tpmod p, k )\mid 0\le k < p\} \tag{$\dagger$}
\end{equation*}
taken under product order. By considering the family of $d$-parallelotopes whose characteristic vertices can be obtained by switching some of the $t_i$'s in $(t_1,t_2,\dots,t_{d-1},p)$ to $p-t_i$'s, we see that, at a loss of constant factor exponential in $d$, we may for $\mathcal{P}$ consider not just the width of \eqref{dag2}, but the maximum of the widths of all posets of the form
\begin{equation*}\label{dag3}
   \{((\pm kt_1)\tpmod p, \dots, (\pm kt_{d-1})\tpmod p,k )\mid 0\le k < p\} \tag{$\ddagger$}
\end{equation*}
taken under product order. Just as in Section \ref{sec3}, we will refer to this collection of posets as the family $S_\mathcal{P}$ (see Definition \ref{def43}).

As we enumerate over all such parallelotopes $\mathcal{P}$, we run the risk of counting a given obstructing hypercube in arbitrarily many such parallelotopes. To avoid this, we employ the same technique used in Section \ref{sec3}. Specifically, we restrict the possible obstructing cubes in each parallelotope to the obstructing hypercubes in the upper half of the parallelotope. In other words, we would only consider elements of \eqref{dag3} corresponding to $\frac{p}{2}< k <p$. By Proposition $\ref{prop32}$ \cmnt{write out the generalization}, the width of this restricted poset differs from that of \eqref{dag3} by at most a factor of two, so we may ignore this range restriction on $k$ and consider the entirety of the set.

Recall that the argument used in three dimensions (see Section $\ref{sec3}$) visualized the elements of the given partially ordered set as points within a cube of side length $p$, from which an antichain could be viewed as the intersection between $[0,p)^3$ and a lattice whose normal vector had coordinates of uniform sign. In generalizing the results from three dimensions into $d>3$ dimensions, we similarly view the elements of our partially ordered sets as points within a lattice. However the similarities end there, for the argument in Section $\ref{sec3}$ does not extend directly to higher dimensions, and so a new approach must be taken.

\subsection{Lattices and the LLL Lattice Basis Reduction Algorithm}
We first recall several important properties of lattices. 


\begin{prop}\label{prop41}
The covolume of a lattice is independent of the choice of basis.
\end{prop}

\begin{prop}\label{prop43} (Hadamard's Inequality)
Let $B=(b_1,\dots,b_n)$ be a basis for a lattice $\mathcal{L}$ in $\mathbb{R}^m$. Then
\[ d(\mathcal{L})\le \prod_{i=1}^n |b_i|,\]
where $d(\mathcal{L})$ is the covolume of $\mathcal{L}$. Equality holds if and only if the elements of $B$ are orthogonal.
\end{prop}

For any basis $B=(b_1,b_2,\dots, b_d)$ of a lattice $\mathcal{L}$, let $B^{*}=(b_1^{*},b_2^{*},\dots,b_d^{*})$ be the result when the Gram-Schmidt orthogonalization procedure is applied to $B$ (without normalizing the vectors $b_i^*$, so the product of their lengths is the covolume of the lattice), and let $\mu_{i,j}=\frac{b_i\cdot b_j^*}{b_j^*\cdot b_j^*}$ be the orthogonal projection coefficient of $b_i$ onto $b_j^*$. We now define an LLL reduced basis \cite{LLL}.
\begin{definition}
The basis $B=(b_1,b_2,\dots, b_d)$ is said to be $\textit{LLL reduced}$ if the following two conditions are met
\begin{enumerate}[label=(\roman*)]
    \item $|\mu_{i,j}| \le \frac{1}{2}$ for all $1\le j < i \le d,$ 
    \item $ |b_i^*+\mu_{i,i-1}b_{i-1}^*| ^2\ge \frac{3}{4}|b_{i-1}^*| ^2$ for all $1\le i \le d.$
\end{enumerate}
\end{definition}

One can think of such a basis as being a good approximation of a short orthogonal basis. In \cite{LLL} it is proved that every lattice has an LLL reduced basis (in fact, an efficient algorithm for finding such a basis is given). If a basis is LLL reduced, Proposition \ref{prop43} can be strengthened to the following proposition, taken from \cite{LLL} [Proposition $1.6$].

\begin{prop}[Lenstra, Lenstra, Lovász \cite{LLL}] \label{prop42}
Let $B=(b_1,\dots,b_d)$ be an LLL reduced basis for a lattice $\mathcal{L}$ in $\mathbb{R}^d$. Then we have
\[d(\mathcal{L})\le \prod_{i=1}^n  |b_i| \le 2^\frac{d(d-1)}{4} d(\mathcal{L}),\]
where $d(\mathcal{\mathcal{L}})$ is covolume of $\mathcal{L}$.
\end{prop}

\subsection{An LLL Reduced Basis of a Familiar Lattice}
Let $p$ be a prime and let $\mathcal{P}$ be the paralleletope with characteristic vertex $(t_1,\cdots,t_{d-1},p)$.

\begin{definition}\label{def42}
Let $\mathcal{L}_\mathcal{P}$ be the lattice whose basis consists of the vector $b_0=(t_1,t_2,\dots,t_{d-1},1)$, and $d-1$ vectors of the form $b_{k}=(0,\dots,0,p,0,\dots,0)$, where the $k$-th such vector has all zeros except for a $p$ in the $k$-th coordinate starting from the left for $1\le k \le d-1$.
\end{definition}
\begin{lemma}\label{lem41}
The covolume of $\mathcal{L}_\mathcal{P}$ is $p^{d-1}$.
\end{lemma}
\begin{proof}
Upon expansion of the determinant of the matrix given by taking the basis given in Definition \ref{def42} as rows
\[ 
\begin{vmatrix}
    t_1 & t_2 & t_3 & \dots & t_{d-1} & 1 \\
    p & 0 & 0 & \dots  & 0 & 0 \\
    \vdots & \vdots & \vdots & \ddots & \vdots & \vdots \\
    0 & 0 & 0 & \dots  & p & 0
\end{vmatrix}
,\]
 all terms equal zero with the exception of one times the determinant of the $(d-1)$ by $(d-1)$ diagonal matrix whose diagonal entries are all $p$, which evaluates to $\pm p^{d-1}$, with sign dependent on the parity of $d$. The lemma then follows by taking absolute value.
\end{proof}

Furthermore, observe that the intersection of $\mathcal{L}_\mathcal{P}$ with $[0,p)^d$ is precisely the multiples of the characteristic vertex  $(t_1,t_2,\dots, t_{d-1},1)$ taken modulo $p$ and as points in $\mathbb{R}^d$. 

\begin{definition}\label{def43}
Let $S_\mathcal{P}$ denote the set of partially ordered sets of the form
\[\big\{ \big( (\pm t_1 \cdot k)\tpmod p, \dots, (\pm t_{d-1} \cdot k)\tpmod p, k \big) \mid 0\le k < p \big \}\] for some fixed choice of signs, and under product order.
\end{definition}

\begin{lemma}\label{lem42}
When the elements of any $S\in S_\mathcal{P}$ are taken as points in $d$ dimensional space, there exists a lattice of points from $S$ whose intersection with $[0,p)^d$ is of size $\Omega (p^{\frac{d-1}{d}})$.
\end{lemma}
\begin{proof}
By Lemma \ref{lem41}, the covolume of $\mathcal{L}_\mathcal{P}$ is $p^{d-1}$. Let $B=(b_1,b_2,\dots,b_d)$ be an LLL reduced basis of $\mathcal{L}_\mathcal{P}$. By Proposition $\ref{prop42}$ we have
\[ p^{d-1} \le \prod_{i=1}^{d} |b_i| \le 2^{\tfrac{d(d-1)}{4}}p^{d-1}. \tag
{\S} \]

Reorder the elements of the basis $B$ such that $| b_i |  \le | b_{i+1} |$ for all $1\le i \le d-1$ - note that this doesn't change the product of the $|b_i|$s. It follows then that
\[
\prod_{i=1}^{d-1}| b_i| \le (2^{\frac{d(d-1)}{4}}p^{d-1})^\frac{d-1}{d}=2^\frac{(d-1)^2}{4}p^{\frac{(d-1)^2}{d}}.
\]
It follows from Proposition \ref{prop43} that the covolume of the fundamental region determined by the first $d-1$ vectors of $B$ is at most $2^\frac{(d-1)^2}{4}p^{\frac{(d-1)^2}{d}}$.



Suppose that $b_i$ contains a coordinate whose absolute value is at least $\frac{p}{2d}$ for $1\le i \le d-1$. It follows that each such $b_i$ has magnitude at least $\frac{p}{2d}$. Referring back to $(\S)$, this would imply that
\[\bigg(\frac{p}{2d}\bigg)^{d-1}\le 2^\frac{(d-1)^2}{4}p^\frac{(d-1)^2}{d} \rightarrow \frac{p}{2d} \le 2^\frac{d-1}{4}p^\frac{d-1}{d}.\]
However $d$ is fixed, and so for sufficiently large $p$, the inequality fails, implying that for sufficiently large $p$, there exists a $b_i$ whose magnitude, and therefore largest  coordinate absolute value wise is at most $\frac{p}{2d}$. As our $b_i$ are sorted by increasing magnitude suppose that for all $1\le i \le k$, $| b_i| \le \frac{p}{2d}$. It follows that
\begin{align}
    \prod_{i=1}^k |b_i| &\le (2d)^{(d-1-k)}2^\frac{(d-1)^2}{4}p^{\tfrac{(d-1)^2}{d}-(d-1-k)}\\
    &=O\big(p^{\tfrac{(d-1)^2}{d}-(d-1-k)}\big)
\end{align}
and so by Proposition \ref{prop43}, we see that the covolume of the fundamental region spanned by the first $k$ basis vectors is $O\big(p^{\tfrac{(d-1)^2}{d}-(d-1-k)}\big)$. 

In a manner analogous to that used in Lemma \ref{lem32}, it can be shown that there exists an element $\ell \in \mathcal{L}_\mathcal{P}$ within the region $\big[\frac{p}{2d},\frac{p(2d-1)}{2d}\big]^d$. 
 It follows that the sublattice $\mathcal{L}_\ell$ of $\mathcal{L}_\mathcal{P}$ obtained by adding $\ell$ to every element of the sublattice spanned by the first $k$ elements of $B$ contains points within the region $(\frac{p}{2d},\frac{p(2d-1)}{2d})^d$. 
 Finally, the number of points in the intersection of $\mathcal{L}_\ell$ with $[0,p)^d$ is
    \[\frac{p^k}{O\big(p^{\tfrac{(d-1)^2}{d}-(d-1-k)}\big)}=\Omega(p^\frac{(d-1)}{d}).\] 
The lemma then follows, as $\mathcal{P}$ can be taken to be any element of $S_\mathcal{P}$.
\end{proof}

\begin{lemma}\label{lem43}
There exists an $S\in S_\mathcal{P}$ for which the width of $S$ is $\Omega (p^{\frac{d-1}{d}})$.
\end{lemma}

\begin{proof}
This proof proceeds analogously to Lemma \ref{lem35}. Recall that $S_\mathcal{P}$ is the set of $2^{d-1}$ posets all differing from each other by the negation of some subset of their coordinates. Furthermore, if $\mathcal{L}_1$ is a lattice such that $\mathcal{L}_1\cap[0,p)^d$ consists only of points within a given poset $S_1\in S_\mathcal{P}$, then $\mathcal{L}_1\cap[0,p)^d$ corresponds to an antichain of $S_1$ if there exists a vector $v$ orthogonal to $\mathcal{L}_1$ such that the coordinates of $v$ are of uniform sign. 

Now observe that just as there exists a lattice $\mathcal{L}_1$ of points from $S_1$, there exist analogous $\mathcal{L}_i$'s for each of the other $2^{d-1}-1$ posets $S_i\in S_\mathcal{P}$ obtained by applying a map taking each element  $\ell\in\mathcal{L}_1$ to $ell'\in \mathcal{L}_i$ with some of the signs reversed.
Crucially, these analogous lattices also have analogous normal vectors, albeit with the signs changed to match those of the chosen $\mathcal{L}_i$. As all $2^{d-1}$ sign combinations are present in $S_\mathcal{P}$, it follows that there exists an element $S_i\in S_\mathcal{P}$ for which the lattice $\mathcal{L}_i$ constructed above has a normal vector all of whose coordinates are of the same sign. The result then follows from the fact that by Lemma \ref{lem42} $|\mathcal{L}_1\cap[0,p)^d|=\Omega (p^{\frac{d-1}{d}})$.
\end{proof}

\begin{remark}
If the elements of $S_\mathcal{P}$ were to behave in the same random manner as that of Theorem \ref{thm22}, their widths would also be $\Omega(p^{\frac{d-1}{d}})$.
\end{remark}

\begin{theorem}\label{thm41.5}
Let $\vec{t} = (t_1, ..., t_{d-1}) \in \{1,...,p-1\}^{d-1}$. Let $S_{\vec{t}}$ denote the set of $2^{d-1}$ partially ordered sets, each of which is of the form
\[
\big\{ \big( (\pm t_1 \cdot k) \tpmod p, \dots, (\pm t_{d-1} \cdot k)\tpmod p, k
\big)\mid 0\le k <p\big\},
\]
under product order, for some fixed choice of signs.
Then for each $\vec{t}$ there exists an element of $S_{\vec{t}}$ whose width is $\Omega (p^{1-\frac{1}{d}})$.
\end{theorem}
\begin{proof}
The result follows directly from Lemma \ref{lem43}.
\end{proof}

\subsection{Conclusion}
We are now ready to state the final theorem of the section which combines previous results to provide a lower bound on the number of simultaneously visible obstructing unit $d$-hypercubes within a $d$-hypercube of sidelength $n$.

\begin{theorem}\label{thm42}
For all positive integers $n>3$, the maximum number of unit $d$-hypercubes visible from an observing unit $d$-hypercube and all within a $d$-hypercube of side length $n$ is $\Omega(n^{d-\frac{1}{d}})$.
\end{theorem}

\begin{proof}
By Bertrand's postulate, we may without loss of generality assume that $n=p$ is a prime. There are $\Theta(p^{d-1})$ families of $2^{d-1}$ paralleletopes. By Theorem \ref{thm41.5} each such family has at least one element, call it $\mathcal{P}$, whose associated partially ordered set has width $\Omega(p^{1-\frac{1}{d}})$. As a result Proposition \ref{prop32}, there is an antichain of this partially ordered set consisting only of points with last coordinate larger than $\tfrac{p-1}{2}$ also of size $\Omega(p^{1-\frac{1}{d}})$. By the bijection between sets of simultaneously visible obstructing cubes along the front lateral edge $E_\mathcal{P}$ of $\mathcal{P}$ and antichains of $\mathcal{P}$'s associated partially ordered set, there must exist a corresponding set of $\Omega(p^{1-\frac{1}{d}})$ simultaneously visible obstructing cubes along $E_\mathcal{P}$ all of which are in the upper half of the paralleletope. By only considering obstructing cubes in the upper half of each paralleletopes, and spacing out the characteristic vertices of the parallelotopes sufficiently, we guarantee that none of these cubes intersect more than one parallelotope at a time. Finally, enumerating over these $\Theta(p^{d-1})$ families of paralleletopes yields a total of $\Omega(p^{d-\frac{1}{d}})$ total visible obstructing unit $d$-hypercubes.
\end{proof}

%% file: Sections/S5.tex
\section{A Toy Upper Bound} \label{sec5}
In the following three sections, we work within the frame of the parallelepiped model we used for the lower bound. Whereas in the lower bound we bounded from below the largest number of cubes that could possibly be seen from within the confines of each parallelepiped, we now bound this value from above. We do so by making use of the same partially ordered sets as used in the lower bound (see Sections \ref{sec3} and \ref{sec4}), this time obtaining an upper bound of the width by constructing chain covers. By Dilworth's Theorem, the size of these chain covers will then serve as upper bounds on the widths of the posets. 


In addition, the following three sections will treat visibility as it applies in the general case of $d>2$ dimensions. We continue with the same notation used in Section \ref{sec4}. 

\subsection{Setup of the Toy Upper Bound}
We supplement the preexisting notation with several additional definitions that will be used heavily in the next two sections.

\begin{definition}\label{def50} For integers $t_i$ and $d > 2$, we denote:
\[\vec{t} := (t_1, t_2, \dots, t_{d-1}, 1).\]
\end{definition}

\begin{definition}\label{def51} For prime $p$, positive integer $d>2$ and $\vec t=(t_1,\dots,t_{d-1},1)\in \mathbb{Z}_p^d$, we define:
\[h_{p}(\vec{t}) := \min_{0 < a < p} \max ((at_1)\tpmod p,\dots,(at_{d-1})\tpmod p, a).\]
\end{definition}

\begin{remark}
The above definition is motivated by the notion of height in projective space. (Refer to \cite{article}.)
\end{remark}

For parallelotope $\mathcal{P}$ with characteristic vertex $(t_1,\cdots,t_{d-1},p)$, we will let $S_{\vec{t}}$ be the partially ordered set associated with $\mathcal{P}$.

\begin{lemma}\label{lem51}
The width of $S_{\vec{t}}$ is at most $d h_p(\vec{t})$.
\end{lemma}

\begin{proof}
We generate a chain cover of $S_{\vec{t}}$ from the first $p-1$ multiples of $\vec{u}$, which we choose to be an element of $S_{\vec{t}}$ whose maximum coordinate is $h_p(\vec{t})$. Writing down all the multiples in the order they appear, we traverse this list from its start. On step one, we create a chain and add the first multiple, namely $\vec{u}$, into it. On step $k$, we examine the tuple $k\cdot \vec{u}$. If each of the coordinates of $k\cdot \vec{u}$ is greater modulo $p$ than its corresponding coordinate in $(k-1) \cdot \vec{u}$, then we append $k\cdot \vec{u}$ onto the end of the current chain. Otherwise, we terminate the current chain and cast it aside, adding $k\cdot \vec{u}$ to a new chain. As $p$ is prime, the multiples of $\vec{u}$ will take on every value in $S_{\vec{t}}$. It follows that after step $p-1$, the collection of chains formed by the process forms a chain cover on $S_{\vec{t}}$ ($0\cdot \vec{u}$ can be appended to the beginning of any of the antichains). Note that an existing chain is completed and a new chain is started at step $k$ if and only if one of the coordinates in the transition from $(k-1)\cdot \vec{u}$ to $k\cdot \vec{u}$ exceeds $p$ and ``loops back'' to a smaller value modulo $p$. The number of steps where this occurs in at least one coordinate is at most the sum of the number of times it occurs in each coordinate, which is equal to $u_1+u_2+\dots + u_d \le d h_p(\vec{t})$. As the size of any chain cover of $S_{\vec{t}}$ is greater than the width of $S_{\vec{t}}$, we are done.  
\end{proof}

As a result of Lemma \ref{lem51}, it suffices to place an upper bound on $h_p(\vec{t})$.

%% file: Sections/S6.tex
\section{The Discrete Fourier Transform}\label{sec6}
 Before we proceed, we take the time to familiarize the reader with several important notions and analytic techniques that will play crucial roles in the calculations of Section \ref{sec7}. In the proof of the upper bound of our reduced visibility problem, we use techniques from Fourier analysis. Below, we provide statements and proofs of the theorems that we will use.

\begin{definition}\label{def61}
Let $e_p(x):=e^{\frac{2\pi{i}x}{p}}$. Consider some function $f:\mathbb{Z}_p^k\longrightarrow \mathbb{C}$. We define the discrete Fourier transform of $f$, denoted as $\hat{f}$ as follows
\[\hat{f}(\vec x):= \sum_{\vec w\in\mathbb{Z}_p^k}^{} e_p(\vec w\cdot \vec x)f(\vec w).\]
\end{definition}

\begin{prop}[Classical]\label{thm61}
For any two functions $f,g:\mathbb{Z}_p^k\longrightarrow \mathbb{C}$, the following holds
\[\sum_{\vec w\in\mathbb{Z}_p^k}^{}f(\vec w)\overline{g(\vec w)}=\frac{1}{p^k}\sum_{\vec x\in\mathbb{Z}_p^k}^{}\hat{f}(\vec x)\overline{\hat{g}(\vec x)}.\]
\end{prop}

Additionally, we introduce the following lemma to be used later.
\begin{lemma}\label{lem61}
For $n\le p$ define ${h}(x)$ to be equal to one if $0\le x < n$ and zero otherwise. Then 
\[\left| {\hat h}(x)\right | \le \min\bigg(\frac{p}{2\left| x \right|},n\bigg). \]
\end{lemma}
\begin{proof}
We first expand the left hand side
\[\left| {\hat h}(x)\right |=\Bigg| \sum_{0 \le w < n} e_p(wx)\Bigg|=\left| \frac{e_p(nx)-1}{e_p(x)-1}\right | = \frac{\left| e_p(nx)-1 \right |}{\left| e_p(x)-1 \right |}\le \frac{2}{\left| e_p(x)-1 \right |}.\]

Now note that $\left| e_p(x)-1 \right |$ is at least $\frac{2}{\pi}$ times the length of the arc subtended by $\frac{2\pi \left| x\right|}{p}$ radians. Thus $\left| e_p(x)-1 \right |\ge\frac{2}{\pi}\cdot\frac{2\pi \left| x\right |}{p}=\frac{4\left| x\right |}{p}$, and so we have that $\frac{2}{\left| e_p(x)-1 \right |}\le \frac{p}{2\left| x \right|}.$ To finish, observe that 
\[\left| {\hat h}(x)\right |=\Bigg| \sum_{0\le w < n} e_p(wx)\Bigg| \le \sum_{0 \le w < n} \left| e_p(wx) \right|=n,\]
by the triangle inequality. It follows that $\left| {\hat h}(x)\right | \le \min\big(\frac{p}{2\left| x \right|},n\big)$ is desired.
\end{proof}

We now generalize Lemma \ref{lem61} to higher dimensions. 
\begin{lemma}[Generalization of Lemma \ref{lem61}]\label{lem62}
For $n\le p$ let $h:\mathbb{Z}_p^d\mapsto \mathbb{C}$ be such that $h(\vec x)=1$ if for every coordinate $x_i$ of $\vec x$, $0\le x_i < n$, and let $h(\vec x)$ be zero otherwise. Then
\[\left| {\hat h}(x_1,\dots, x_d)\right |\le \prod_{k=1}^{d} \min\bigg(\frac{p}{2\left| x_k \right|},n\bigg).\]
\end{lemma}
\begin{proof}
It suffices to note that 
\[\left| {\hat h}(x_1,\dots, x_d)\right |=\Bigg| \sum_{\vec w\in [0,n)^d} e_p(\vec w\cdot \vec x)\Bigg |=\prod_{k=1}^{d}\Bigg|\sum_{0\le w <n} e_p(wx_k)\Bigg|\le\prod_{k=1}^{d} \min\bigg(\frac{p}{2\left| x_k \right|},n\bigg),\]
where the last inequality is Lemma \ref{lem61} applied to each coordinate of the $x_i$.
\end{proof}

%% file: Sections/S7.tex
\section{Proof of Toy Upper Bound} \label{sec7}
In this section we compute a bound on $h_p(\vec{t})$ (see Definition \ref{def51} for the definition of $h_p(\vec{t})$). For brevity we will refer to this value as $h_p$.
\begin{definition}\label{def71}
Let $f: \mathbb{Z}_p^d\mapsto \{0,1\}$ be such that $f(\vec w)$ equals $1$ if $\vec w$ is a scalar multiple of $\vec t$ modulo $p$ (see Definition \ref{def50}) and zero otherwise.
\end{definition}



\begin{lemma}\label{lem71}
For any $\vec x\in \mathbb{Z}_p^d$, the following holds:
\[\hat{f}(\vec x) =\begin{cases} p & \text{if } \vec x \cdot \vec t=0 \\ 0 & \text{else.} \end{cases} \]
\end{lemma}

\begin{proof} This follows immediately from the formula
\[
\hat{f}(\vec x)= \sum_{\vec{w}} e_p(\vec w\cdot \vec x)f(\vec w)=\sum_{k=0}^{p-1} e_p(k\vec{t}\cdot \vec{x}).\qedhere
\]

\end{proof}
\begin{definition}\label{def72}
We say that $\textbf{1}_S:\mathbb{Z}_p^d\mapsto \{0,1\}$ is the \textit{indicator function} for the set $S$ if for all $\vec w\in \mathbb{Z}_p^d$, $\textbf{1}_S(\vec w)$ equals one if $w \in S$, and zero otherwise.
\end{definition}
\begin{definition}\label{def73}
Let $g_{(d,k)}:\mathbb{Z}_p^d\mapsto \mathbb{R}$ be the indicator function $\textbf{1}_{[0,\frac{h_p}{k})^d}$ convolved with itself $k$ times:
\[
g_{(d,k)}=\textbf{1}_{[0,\frac{h_p}{k})^d} \ast \dots \ast \textbf{1}_{[0,\frac{h_p}{k})^d}.
\]
\end{definition}
\begin{lemma}\label{lem72}
$| \hat g_{(d,k)} (\vec 0) |= \ceil{\frac{h_p}{k}}^{kd}.$
\end{lemma}
\begin{proof}
Upon expansion, we have
\[
|\hat g_{(d,k)} (\vec 0) | = \sum_{\vec w \in \mathbb{Z}_p^d} g_{(d,k)}(\vec w)=\sum_{\vec w \in \mathbb{Z}_p^d}\sum_{\vec v_1+\dots +\vec v_k=\vec w}\prod_{i=1}^{k}\textbf{1}_{[0,\frac{h_p}{k})^d}(v_i)=\sum_{v_1,\dots, v_k\in [0,\tfrac{h_p}{k})^d} 1 = \Big\lceil\frac{h_p}{k}\Big\rceil^{kd}.\qedhere
\]
\end{proof}

More generally, we have the following bound on $|\hat{g}|$.

\begin{lemma}\label{lem73}
Let $\vec x=(x_1,\dots, x_d)$. Then
\[
|\hat g_{(d,k)}(\vec x) | \le \prod_{i=1}^{d} \min\bigg ( \frac{p}{2|x_i|} ,\Big\lceil\frac{h_p}{k}\Big\rceil \bigg)^k \le \bigg(\frac{p}{2\max |x_i |} \Big\lceil\frac{h_p}{k}\Big\rceil^{d-1} \bigg )^k.
\]
\end{lemma}
\begin{proof}
We have
\[
\hat g_{(d,k)}(\vec x) =\sum_{\vec w \in \mathbb{Z}_p^d} e_p(\vec w \cdot \vec x)\sum_{\vec v_1+\dots +\vec v_k=\vec w} \prod_{i=1}^k \textbf{1}_{[0,\tfrac{h_p}{k})^d}(\vec v_i) = \sum_{\vec v_1,\dots, \vec v_k \in [0,\tfrac{h_p}{k})^d} e_p(\vec x \cdot \Sigma \vec v_i).
\]
Rewriting $e_p(\vec x \cdot \Sigma \vec v_i)=\prod_{i=1}^k e_p(\vec x \cdot \vec v_i)$, this reduces to
\[\sum_{\vec v_1,\dots, \vec v_k \in [0,\tfrac{h_p}{k})^d} \prod_{i=1}^k e_p(\vec x \cdot \vec v_i)=\bigg(\sum_{\vec v\in [0,\tfrac{h_p}{k})^d} e_p(\vec x\cdot \vec v) \bigg)^k=\bigg(\hat{\textbf{1}}_{[0,\tfrac{h_p}{k})^d}(\vec x)\bigg)^k.\]
The left inequality now follows from Lemma \ref{lem62} applied to $\textbf{1}_{[0,\frac{h_p}{k})^d}(\vec x)$.
\end{proof}

\begin{definition}\label{dualhp} For $\vec{t} \in \mathbb{Z}_p^d$, we define the \emph{dual height} $h_p^*(\vec{t})$ by
\[
h_p^*(\vec{t}) = \displaystyle{\min_{\substack{\vec \alpha\cdot \vec{t} \equiv 0 \\ \alpha \not \equiv \vec 0}} \max |\alpha_i|}.
\]
\end{definition}

\begin{lemma}\label{lemhp} For all $\vec{t} \in \mathbb{Z}_p^d$, we have
\[
h_p(\vec{t}) \le e(d-1)\frac{p\ceil{\log p}}{2h_p^*(\vec{t})} \ll \frac{p\log p}{h_p^*(\vec{t})},
\]
where $e$ is the base of the natural logarithm.
\end{lemma}
\begin{proof}
Observe that for a fixed $\vec{t}$, if $f$ is defined as in Definition \ref{def71}, we have
\[
\sum_{\vec w \in \mathbb{Z}_p^d} f(\vec w) \overline {g_{(d,k)}(\vec w)}=\sum_{\vec w \in \mathbb{Z}_p^d} f(\vec w) g_{(d,k)}(\vec w)=f(\vec{0})g_{(d,k)}(\vec{0})=1,
\]
where the middle equality used the fact that by the definition of $h_p(\vec{t})$, the support of $g_{(d,k)}$ only intersects the support of $f$ at $\vec{0}$. However by Proposition \ref{thm61}, we also have 
\[
p^d\sum_{\vec w \in \mathbb{Z}_p^d} f(\vec w) \overline {g_{(d,k)}(\vec w)}=\sum_{\vec \alpha \in \mathbb{Z}_p^d} \hat f(\vec \alpha) \overline {\hat g_{(d,k)}(\vec \alpha)}.
\]
Note that $\hat f$ is only supported on $\alpha$ if $\vec \alpha\cdot \vec{t} \equiv 0 \tpmod{p}$, hence
\[
p^d=\sum_{\vec \alpha \in \mathbb{Z}_p^d} \hat f(\vec \alpha) \overline {\hat g_{(d,k)}(\vec \alpha)} = \sum_{\vec \alpha\cdot \vec{t} \equiv 0} \hat f(\vec \alpha) \overline {\hat g_{(d,k)}(\vec \alpha)} =\sum_{\vec \alpha\cdot \vec{t} \equiv 0} p\overline {\hat g_{(d,k)}(\vec \alpha)},
\]
where the last equality is due to Lemma \ref{lem71}, so
\[p^{d-1}=\sum_{\vec \alpha\cdot \vec{t} \equiv 0} \hat g_{(d,k)}(\vec \alpha).
\]
Dividing into cases based on whether $\vec{\alpha} = 0$, we have 
\[
\hat g_{(d,k)} (\vec 0) + \sum_{\substack{\vec \alpha\cdot \vec{x} \equiv 0 \\ \alpha \not \equiv \vec 0}} {\hat g_{(d,k)}(\vec \alpha)}=p^{d-1},
\]
so by the triangle inequality we have
\begin{equation*}\label{dag1}
\hat g_{(d,k)} (\vec 0) \le p^{d-1} +\sum_{\substack{\vec \alpha\cdot \vec{t} \equiv 0 \\ \alpha \not \equiv \vec 0}}|\hat g_{(d,k)}(\vec \alpha)| \tag{$\dagger$}.
\end{equation*}
For each $\vec{\alpha}$, let $\vec{\alpha}=(\alpha_1,\dots,\alpha_d)$. Plugging the results of lemmas \ref{lem72} and \ref{lem73} into \eqref{dag1}, we have 
\[
\Big\lceil\frac{h_p}{k}\Big\rceil^{dk} \le p^{d-1} +\sum_{\substack{\vec \alpha\cdot \vec{t} \equiv 0 \\ \alpha \not \equiv \vec 0}}|\hat g_{(d,k)}(\vec \alpha)|\le p^{d-1} +\sum_{\substack{\vec \alpha\cdot \vec{t} \equiv 0 \\ \alpha \not \equiv \vec 0}} \bigg(\frac{p}{2\max |\alpha_i |} \Big\lceil \frac{h_p}{k}\Big\rceil^{d-1} \bigg )^k.
\]
Dividing both sides by $\ceil{\frac{h_p}{k}}^{(d-1)k}$ yields 
\begin{align*}
\Big\lceil\frac{h_p}{k} \Big\rceil^k & \le \frac{p^{d-1}}{\ceil{\frac{h_p}{k}}^{(d-1)k}}+\sum_{\substack{\vec \alpha\cdot \vec{t} \equiv 0  \\ \alpha \not \equiv \vec 0}} \bigg(\frac{p}{2\max |\alpha_i |}\bigg)^k \\
& \le \frac{p^{d-1}}{\ceil{\frac{h_p}{k}}^{(d-1)k}} + (p^{d-1}-1)\max_{\substack{\vec \alpha\cdot \vec{t} \equiv 0 \\ \alpha \not \equiv \vec 0}} \bigg(\frac{p}{2\max |\alpha_i | } \bigg) ^k.
\end{align*}
We now let $h_p^*(\vec{t}) = \displaystyle{\min_{\substack{\vec \alpha\cdot \vec{t} \equiv 0 \\ \alpha \not \equiv \vec 0}} \max |\alpha_i|}$. Then
\begin{align*}
\Big\lceil\frac{h_p}{k} \Big\rceil^k &\leq \frac{p^{d-1}}{\ceil{\frac{h_p}{k}}^{(d-1)k}} + (p^{d-1}-1)\Big(\frac{p}{2h_p^*} \Big)^k \\
&\leq p^{d-1} \max \bigg( \frac{p^{d-1}}{\ceil{\frac{h_p}{k}}^{(d-1)k}}, \Big( \frac{p}{2h_p^*} \Big)^k \bigg).
\end{align*}
We consider two cases. 

{\bf Case 1: } First, suppose that
\[
\max \bigg( \frac{p^{d-1}}{\ceil{\frac{h_p}{k}}^{(d-1)k}}, \Big( \frac{p}{2h_p^*} \Big)^k \bigg) = \frac{p^{d-1}}{\ceil{\frac{h_p}{k}}^{(d-1)k}}.
\]
Rearranging yields
\[
\Big\lceil\frac{h_p}{k}\Big\rceil^k \leq p^{d-1} \max \bigg( \frac{p^{d-1}}{\ceil{\frac{h_p}{k}}^{(d-1)k}}, \Big( \frac{p}{2h_p^*} \Big)^k \bigg) =  \frac{p^{2d-2}}{\ceil{\frac{h_p}{k}}^{(d-1)k}},
\]
so using $\frac{h_p}{k} \le \ceil{\frac{h_p}{k}}$, we get
\[
\bigg( \frac{h_p}{k} \bigg)^{dk} \leq p^{2d-2},
\]
so
\[
h_p \le k\cdot p^{\tfrac{2d-2}{dk}},
\]
and from the trivial inequality $2h_p^* \le p$, we get
\[
h_p\cdot 2h_p^* \le k\cdot p\cdot p^{\tfrac{2d-2}{dk}} \le k\cdot p\cdot p^{\tfrac{d-1}{k}},
\]
This concludes the first case. $\clubsuit$

{\bf Case 2: } We suppose that
\[
\max \bigg( \frac{p^{d-1}}{\ceil{\frac{h_p}{k}}^{(d-1)k}}, \Big( \frac{p}{2h_p^*} \Big)^k \bigg) =  \Big( \frac{p}{2h_p^*} \Big)^k.
\]
Then rearranging the inequality gives:
\[
\Big\lceil \frac{h_p}{k} \Big\rceil^k \leq p^{d-1} \max \bigg( \frac{p^{d-1}}{\ceil{\frac{h_p}{k}}^{(d-1)k}}, \Big( \frac{p}{2h_p^*} \Big)^k \bigg) \leq p^{d-1} \Big( \frac{p}{2h_p^*} \Big)^{k},
\]
so
\[
\Big\lceil \frac{h_p}{k} \Big\rceil^k \leq p^{d-1} \cdot \Big( \frac{p}{2h_p^*} \Big)^{k},
\]
which simplifies to
\[
h_p \cdot 2h_p^* \leq k\cdot p\cdot p^{\tfrac{d-1}{k}}.
\]
This concludes the second case. $\clubsuit$

Now we set $k = \ceil{(d-1)\log p}$ yielding
\[
h_p \cdot 2h_p^* \leq \ceil{(d-1)\log p} \cdot p \cdot p^{\tfrac{d-1}{(d-1)\log p}} \le e\cdot (d-1) \cdot p\ceil{\log p}.\qedhere
\]
\end{proof}

\begin{theorem}\label{thm74}
The average value of $h_p(\vec{t})$ as $\vec{t}$ varies is bounded by $O(p^{\tfrac{d-1}{d}} \log p)$, that is,
\[
\mathbb{E}_{\vec{t}}[h_p(\vec{t})] \ll p^{\tfrac{d-1}{d}} \log p.
\]
\end{theorem}
\begin{proof}
We introduce the parameter $x$, a bound on the size of $h_p^*$. Observe that
\[ h_p^* \leq x \implies \frac{p \log p}{h_p^*} \geq \frac{p \log p}{x},\]
while
\[h_p^* \geq x \implies \frac{p \log p}{h_p^*} \leq \frac{p \log p}{x}.\]
Let $\mathbb{E}[h_p]$ denote the expected value of $h_p(\vec{t})$ and let $\mathbb{E}\big[\textbf{1}_{h_p^*<x} \cdot \frac{p\log p}{h_p^*}\big]$ denote the expected value of $\frac{p\log p}{h_p^*}$ if $h_p^* < x$. Observe that for each $\vec{\alpha} \not\equiv 0$, there are $p^{d-1}$ tuples $\vec t$ for which $\vec t \cdot \vec \alpha \equiv 0$. It follows that
    \begin{align*}
        \mathbb{E}[h_p] &= \mathbb{E}\Big[\textbf{1}_{h_p^*\ge x} \cdot  \frac{p \log p}{x}\Big] +  \mathbb{E}\Big[\textbf{1}_{h_p^*<x} \cdot  \frac{p \log p}{h_p^*}\Big] \\
        &\ll \frac{p \log p}{x} + \frac{1}{p^d} \sum_{x \geq \alpha_i \geq -x} \frac{p \log p}{\max |\alpha_i|}\cdot p^{d-1} \\
        &\ll \frac{p \log p}{x} + \sum_{x \geq \alpha_i \geq 0} \frac{\log p}{\max \alpha_i}.
    \end{align*}

Assuming that the coordinates of $\alpha$ are strictly descending, in doing so losing at most a constant factor, we have
\begin{align*}
\mathbb{E}[h_p] &\ll \frac{p \log p}{x} + \sum_{x \geq \alpha_1 \geq \alpha_2 \dots \geq \alpha_d \geq 0} \frac{\log p}{\alpha_1} \\
&\ll \frac{p \log p}{x} + \log p \sum_{\alpha_1 = 1}^{x} \frac{1}{\alpha_1}\cdot \alpha_1^{d-1} \\
&\ll \frac{p \log p}{x} + x^{d-1}\log p.
\end{align*}
 \cmnt{prev follows from hockey stick identity}
We may now let $x$ be $p^\frac{1}{d}$, yielding $\mathbb{E}[h_p] \ll p^{\frac{d-1}{d}}\log p$.
\end{proof}

We are now ready to state the main result of this section.
\begin{theorem}\label{thm71}
The largest number of cubes visible in the $d$ dimensional toy upper bound visibility environment is $O(p^{d-\tfrac{1}{d}}\log p)$. 
\end{theorem}
\begin{proof}
By Lemma \ref{thm74}, the largest number of visible obstructions within the average $d$-parallelotope is $O(p^{\tfrac{d-1}{d}} \log p)$. Multiplying by $p^{d-1}$ such parallelotopes yields the desired $O(p^{d-\tfrac{1}{d}} \log p)$.
\end{proof} 

%% file: Sections/S8.tex
\section{Upper Bound on Shallow Sight Visibility} \label{sec8}
In this section we prove our final result. While we aren't able to solve the original problem, we are able to provide a bound on a weaker problem in which light passes through any surface that it hits at a shallow enough angle. 

Often in this section it will be inconvenient to repeatedly write out ``$a \text{ }(\text{mod }p)$'' to denote the residue class of $a$ modulo $p$. Instead, we use the following notation.
 
 \begin{definition} For integers $a$ and $p$, we define $a\%p$ to be the least nonnegative remainder of $a$ modulo $p$, so
 \[a\%p \equiv a \text{ }(\text{mod }p) \text{ and } 0 \leq a\%p \leq p-1.\]
 \end{definition}
 
\begin{definition} \label{def81}
We say that a $d$-dimensional hypercube blocks a ray of light \emph{at the angle} $\theta$ if the ray of light intersects some $d-1$-dimensional facet of the hypercube at an angle at most $90^{\circ}-\theta$ away from the normal vector to that facet.
\end{definition}

Generally we will take $\theta = 45^{\circ}$, but our arguments extend to any fixed $\theta > 0$ at the cost of a constant factor in the bound. In two dimensions, it wasn't necessary to consider this weakening of the problem because of the following easy geometric fact.

\begin{figure}
    \centering
    \includegraphics[scale=0.45]{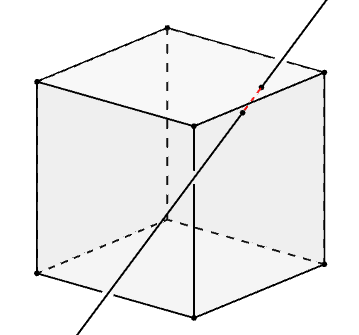}
    \caption{Ray of light passing through a three dimensional cube at a shallow angle.}
    \label{fig:my_label}
\end{figure}

\begin{prop} If $d = 2$, then a square blocks a ray of light at the angle $45^{\circ}$ if and only if it blocks the ray of light.
\end{prop}

With this shallow-angle setup, we can divide the problem into separate bounds for each type of facet, and focus only on $d-1$-dimensional facets which are constant in the last coordinate, where the last coordinate is larger than any of the other coordinates for every point in the facet. Additionally, within each parallelotope we will restrict attention to those facets which intersect a particular edge of that parallelotope - without loss of generality, we consider the edge which passes through the origin $(0,...,0)$.

A second difficulty we face is that unlike the toy problem considered in the previous section, we have to consider light rays going in any possible direction, not just those going in the direction parallel to the long edge of the current parallelotope. The issue is that if one bounds the number of visible obstructions in each parallelotope naively, then the best upper bounds we can hope to prove are much too large, since obstructions which are close to the origin are counted very many times over (since they intersect many parallelotopes). This difficulty already came up in the two-dimensional setting, where the solution was to restrict attention to obstructions that occur in a given parallelotope ``for the first time'', so that each visible obstruction is only counted once. More precisely, we wish to only count an obstruction within the parallelotope which has the largest possible intersection with the obstruction (so that it blocks as many other potential obstructions within that parallelotope as possible).



\begin{definition} We say that a $d-1$-dimensional facet $\mathcal{F}$ (with constant last coordinate) of a unit $d$-dimensional hypercube (axis aligned, with integer coordinates) is a \emph{primitive obstruction} of the parallelotope $\mathcal{P}$ with characteristic vertex $(t_1, ..., t_{d-1}, p)$ if $\mathcal{F}$ intersects the edge connecting the origin to the characteristic vertex of $\mathcal{P}$ and if for every $1 \le i \le d-1$, the line connecting the origin to the vertex $(t_1, ..., t_{i-1}, t_i-1, t_{i+1}, ..., d_{d-1}, p)$ does \emph{not} intersect the facet $\mathcal{F}$.
\end{definition}

\begin{prop} A $d-1$-dimensional facet $\mathcal{F}$ is a primitive obstruction of the parallelotope $\mathcal{P}$ iff it intersects both the line connecting the origin to the characteristic vertex and the line connecting the point $(1,1,...,1,0)$ to the characteristic vertex.
\end{prop}

\begin{prop} Suppose that the facet $\mathcal{F}$ intersects the edge connecting the origin to the characteristic vertex of a parallelotope $\mathcal{P}$ and that some point $p$ of $\mathcal{F} \cap \mathcal{P}$ is not obstructed by any other obstructing facets which intersect that edge of $\mathcal{P}$. Then there is a unique parallelotope $\mathcal{P}'$ such that $\mathcal{F}$ is a primitive obstruction of $\mathcal{P}'$, and the point $p$ will be contained in $\mathcal{F} \cap \mathcal{P}'$ and will not be obstructed by any other obstructing facets that intersect the edge connecting the origin to the characteristic vertex of $\mathcal{P}'$.
\end{prop}

Take the edge $E_{\mathcal{P}}$ connecting the origin and the characteristic vertex of $\mathcal{P}$. We only consider the set of obstructions that intersect $E_{\mathcal{P}}$. The first crucial observation here is that if a facet $\mathcal{F}$ intersects $E_{\mathcal{P}}$, and $V$ is its vertex inside $\mathcal{P}$, then if $V$ is visible from a point of the unit $(d-1)$-hyperface formed by projecting the $d$-hypercube at the origin onto its bottom hyperface, $V$ is visible from the coordinate $(1, 1, \dots, 1, 0)$. Thus, in considering whether an obstruction is visible or not, we only need to consider its vertex $V$ inside $\mathcal{P}$. Note that the coordinates of the vertices in $\mathcal{P}$ belonging to obstructing cubes that intersect $E$ are of the form
\[\Bigg( \bigg\lceil \frac{t_1 \cdot a}{p} \bigg\rceil, \bigg\lceil \frac{t_2 \cdot a}{p} \bigg\rceil, \dots, \bigg\lceil \frac{t_{d-1} \cdot a}{p} \bigg\rceil, a \Bigg).\]
This brings us to the following observation. Consider two obstructions intersecting $E_{\mathcal{P}}$ whose vertices inside $\mathcal{P}$ are $X_1$ and $X_2$, with $X_1$ having a smaller $x_d$-coordinate than $X_2$. If the $x_i$-coordinate slope connecting $X_1$ to $(1, 1, \dots, 1, 0)$ is greater than the $x_i$-coordinate slope connecting $X_2$ to $(1, 1, \dots, 1, 0)$ for every $1 \leq i \leq d-1$, then $X_1$ obstructs $X_2$.

As a result, the upper bound on bottom hyperfaces reduces to computing the width of the following partially ordered set, where the elements are taken under the order $X_1$ $\leq$ $X_2$ if every coordinate in $X_1$ is greater than every coordinate in $X_2$:
\[ \label{dag4}
\Bigg\{ \bigg( \frac{ \lceil \frac{t_1 \cdot a}{p} \rceil -1}{a}, \frac{ \lceil \frac{t_2 \cdot a}{p} \rceil -1}{a}, \dots, \frac{ \lceil \frac{t_{d-1} \cdot a}{p} \rceil -1}{a}, p-a \bigg) \Bigg\} \tag{$\ddagger$}.
\]





\begin{lemma} The partially ordered set from $\eqref{dag4}$ when restricted to primitive obstructions is equivalent to the following poset:
\[ \Bigg\{ \Big(  \frac{ (t_1 \cdot k) \%p }{k}, \frac{ (t_2 \cdot k) \% p }{k}, \dots, \frac{ (t_{d-1} \cdot k) \% p }{k}, k  \Big)   \big | (t_i \cdot k)\% p < k \Bigg\}.\]
\end{lemma}
\begin{proof}
For $1 \leq i \leq d-1$, we have that $\lfloor \tfrac{t_i \cdot a}{p}\rfloor -1 > a \cdot \tfrac{t_i-1}{p}$, which reduces to $(t_i \cdot a) \% p < a$. The poset can then be written as
\[ \Bigg\{ \Big(  \frac{ (t_1 \cdot k) \%p }{k}, \frac{ (t_2 \cdot k) \% p }{k}, \dots, \frac{ (t_{d-1} \cdot k) \% p }{k}, k  \Big)   \mid (t_i \cdot k)\% p < k \Bigg\},\]
where we are inverting the order of the prior poset. That is, $X_1 \leq X_2$ if every coordinate in $X_1$ is less than every coordinate in $X_2$.
\end{proof}

In order to find an upper bound on the width of this poset, we will construct a chain cover of this poset.

To do so, consider vectors of the form:
$$ \bigg( \frac{(t_1\cdot l) \% p}{l}, \frac{(t_2 \cdot l) \%p}{l}, \dots, \frac{(t_{d-1}\cdot l) \% p}{l}, l \bigg),$$
that satisfy the conditions $(t_i\cdot l)\%p \ge l$ for $1 \leq i \leq d-1$. Supposing $l$ satisfies this condition, then for any $a$ such that $(t_i \cdot a)\% p < a$ for all $i$ and such that $(t_i\cdot a)\%p + (t_i\cdot l)\%p < p$ for all $i$, we have
\begin{align*}
\bigg( \frac{(t_1\cdot a) \%p}{a},\dots, \frac{(t_{d-1}\cdot a)\%p}{a} \bigg)& < \bigg( \frac{(t_1\cdot a)\%p+(t_1\cdot l)\%p}{a+l}, \dots, \frac{(t_{d-1}\cdot a)\%p+(t_{d-1}\cdot l)\%p}{a+l} \bigg)\\
 &=\bigg( \frac{(t_1(a+l))\%p}{a+l}, \dots, \frac{(t_{d-1}(a+l))\%p}{a+l} \bigg).
\end{align*}
\cmnt{ in case i screwed things up accidentally
$$\bigg( \frac{i\cdot k \% p}{k}, \frac{j \cdot k}{k} \bigg) + \bigg( \frac{i\cdot l \% p}{l}, \frac{j \cdot l}{l}, l \bigg)$$ 
$$ > \bigg( \frac{(i\cdot k)\%p + (i\cdot l)\%p}{l+k}, \frac{(j\cdot k)\%p + (j\cdot l)\%p}{l+k} \bigg)$$ 
$$> \bigg( \frac{i(l+k)\%p}{l+k}, \frac{j(l+k)\%p}{l+k}\ \bigg).$$ 
}

The condition $(t_i\cdot a)\%p + (t_i\cdot l)\%p < p$ fails for at most $\sum_{i=1}^{d-1} (t_i \cdot l)\%p$ values of $a$. Thus, we can partition the poset into at most $l+\sum_{i=1}^{d-1} (t_i \cdot l)\%p$ chains. It now suffices to prove that for each $(t_1, ..., t_{d-1})$ we can find an $l$ such that $(t_i\cdot l)\% p \ge l$ for all $i$, with the average value of $l+\sum_{i=1}^{d-1} (t_i \cdot l)\%p$ (averaged over all choices of $\vec{t}$) bounded by $O(p^{\tfrac{d-1}{d}}\log p)$.

To find a convenient value of $l$, we choose $l$ such that the maximal coordinate of the vector $(((t_1-1)l)\%p, ((t_2-1)l)\%p, ..., ((t_{d-1}-1)l)\%p, l)$ is at most $h_p(t_1-1, ..., t_{d-1}-1,1)$ (unless $h_p(...) > \tfrac{p}{2}$, in which case we just take $l = 1$). Then for each $i$, we have $(t_i\cdot l)\% p = ((t_1-1)l)\%p + l \ge l$, and
\[
l+\sum_{i=1}^{d-1} (t_i \cdot l)\%p = dl + \sum_{i=1}^{d-1} ((t_i-1)l)\%p \le (2d-1)h_p(t_1-1, ..., t_{d-1}-1,1).
\]
If $h_p(...) > \tfrac{p}{2}$, then we instead have the bound $l+\sum_{i=1}^{d-1} (t_i \cdot l) \% p < (d-1)p < (2d-1)h_p(...)$.




Which brings us to the main result of the upper-bound portion of this paper:

\begin{theorem}\label{thm81}
The largest number of visible obstructions, in the setting where light fails to interact with any obstruction that does not block it at the angle $45^\circ$, is at most $O(n^{d-\frac{1}{d}}\log n)$.
\end{theorem}
\begin{proof}
Summing over the facets counted in each $\mathcal{P}$, whose characteristic vertex's largest coordinate is the $x_d$ coordinate, gives us a count of at most
$O(p^{d-1}) \cdot O(p^{\tfrac{d-1}{d}}\log p)$ = $O(n^{d-\tfrac{1}{d}}\log n)$ unit hypercubes cubes that can be visible from the observer's hypercube at the origin among this set of $\frac{1}{d} p^d$ obstructing cubes. Therefore, the upper bound on the bottom hyperfaces of the above set of obstructions is $O(n^{d-\frac{1}{d}}\log n)$.
\end{proof}

Essentially the same argument works for any angle $\theta > 0$ replacing $45^\circ$, so long as we sum over all parallelotopes $\mathcal{P}$ with characteristic vertices $(t_1, ..., t_{d-1},p)$ satisfying $t_i \le \frac{p}{\tan(\theta)}$ for all $i$. Note that many of these parallelotopes leave the $p\times \cdots \times p$ hypercube which contains our potential obstructing cubes, and the bound degrades by a factor of $\frac{1}{\theta}$ as $\theta \rightarrow 0$.

%% file: Sections/Future_Work.tex
\section{Future Work}

In this paper, a $\Omega(n^{d-\tfrac{1}{d}})$ lower bound on $f_d(n)$ was presented, as well as $O(n^{d-\tfrac{1}{d}}\log n)$ upper bounds in two environments, the first of which saw sight lines being restricted to only those parallel to the edges of the paralleletope in consideration, and second of which saw unit $d$-hypercubes replaced with their bottom hyperfaces (via a shallow light reformulation). Both of these visibility environments simplified the task of working with an upper bound on visibility, but each did so in a manner compromising their ability to make a statement about the original visibility problem. It remains to find an upper bound respective of all the constraints of the original problem. Furthermore, it still remains to determine whether or not the Toy Upper Bound (see Section \ref{sec7}) can be tightened to $O(n^{d-\tfrac{1}{d}})$.


%% file: Sections/Acknowledgements.tex
\section{Acknowledgements}
We would like to thank Professor Pavel Etingof, Dr. Slava Gerovitch, Dr. Tanya Khovanova, and the MIT PRIMES program for allowing us to conduct this project. We also like to thank Dr. Tanya Khovanova, Dr. Svetlana Makarova and Dr. Claude Eicher for their invaluable feedback in the writing of this paper.